\documentstyle[12pt]{article}
\date{}

\newcommand{\be}{\begin{equation}}
      \newcommand{\ee}{\end{equation}}
      \newcommand{\ba}{\begin{eqnarray}}
       \newcommand{\ea}{\end{eqnarray}}
\newcommand{\ban}{\begin{eqnarray*}}
       \newcommand{\ean}{\end{eqnarray*}}

 \newcommand{\qed}{\hspace*{\fill}\rule{3mm}{3mm}\quad}
 \newcommand{\Pf}{\noindent {\em Proof.} }

\newcommand{\sect}[1]{\section{#1} \setcounter{equation}{0}}

\newtheorem{lem}{Lemma}[section]

\begin{document}
\newtheorem{defn}[lem]{Definition}
 \newtheorem{theo}[lem]{Theorem}
 \newtheorem{prop}[lem]{Proposition}
 \newtheorem{rk}[lem]{Remark}
 \newtheorem{ex}[lem]{Example}
 \newtheorem{note}[lem]{Note}
 \newtheorem{conj}[lem]{Conjecture}
 \newtheorem{cor}[lem]{Corollary}

\title{On the $l$-Function and the Reduced Volume of Perelman II
}
\author{Rugang Ye\\Department of Mathematics\\University of California,
Santa Barbara} 
\maketitle


\sect{Introduction}

This paper is a sequel of [Y2]. In [Y2], a number of  geometric and 
analytic properties of the $l$-function and the reduced volume of Perelman 
were derived. In this sequel, we present a major application of the 
$l$-function and the reduced volume, namely their application to the analysis 
of the asymptotical limits of $\kappa$-solutions. Our focus is to present 
a complete and detailed proof of Proposition 11.2 in [P1], which 
provides asymptotical convergence of $\kappa$-solutions and identifies the
asymptotical
limits to be nonflat gradient shrinking solitons. A main foundation for the proof is a study of the limit $l$-function and the asymptotical 
reduced volume (see Definition 1) on the asymptotical limits of $\kappa$-solutions.
This is an extension and further development of the theory presented in [Y2].

Recall that $\kappa$-solutions are $\kappa$-noncollapsed nonflat ancient solutions of the Ricci flow
with bounded and nonnegative curvature operator. 
In dimension three,
rescaling limits of a solution of the Ricci flow near a blow-up singularity are 
$\kappa$-solutions. In higher dimensions they also arise as the same kind of rescaling limits 
under suitable assumptions on the curvature.
In other words, in various situations, blow-up singularities of the Ricci flow
are modelled on $\kappa$-solutions.  For this reason, it is important to 
understand the structures of $\kappa$-solutions. To this purpose, one 
performs blow-down rescalings of a given $\kappa$-solution, extracts 
smooth limits via the process, and identifies them. This is precisely the content 
of [Proposition 11.2, P1].  This proposition plays an important role in 
the analysis of blow-up singularities given in Perelman's papers [P1]
and [P2]. 

Now we give a more detailed account of the content of this paper.
In Section 2, 
we first derive an upper bound for $\Delta l$ (the Laplacian of the $l$-function) in the weak sense.
Then we derive an analytic lemma concerning strong convergence of 
Sobolev functions. These two tools will be needed for proving
[Proposition 11.2, P1], as will be explained below.

Section 3 and Section 4 are devoted to the proof of [Proposition 11.2, P1]. 
The proof follows the 
sketch of proof in [P1] for [Proposition 11.2, P1] and incorporates the ideas in Sections 9 and 
10 in [P1] regarding how to idenfitiy gradient shrinking solitons. 
The statements of the theorem and the proof of the convergence part are given in Section 3, while the 
main part of the identification of the 
asymptotical limits is presented in Section 4.   The proof of the convergence part is 
based on the  estimates in [Theorem 2.16, Y2] (in particular the curvature estmate [(2.48), Y2] which is derived from [(7.16), P1])
and the upper bound for 
$l$ in [Lemma 3.2, Y2].  

The identification of  
the asymptotical limits can be compared with [Theorem 4.9, Y2], which 
says that a solution of the backward Ricci flow must be a
gradient shrinking soliton, provided that the values of its 
reduced volume at two different times are equal. Recall that, basically, 
the said equality forces the weakly formulated differential inequality 
[(2.68 ), Y2]  (i.e. [(7.13), P1]) to become an equality, which then leads to the smoothness of 
the $l$-function and implies that 
the differential inequality [(2.69), Y2](i.e. [(7.14), P1]) also  becomes an equality. These two equalities 
combined with  
[Lemma 4.8, Y2] (i.e. [Proposition 9.1, P1]) then lead to the desired gradient shrinking soliton 
conclusion. We apply this strategy to the current situation.  First, in Section 3, we apply 
Perelman's lower bound for $l$  as presented in [Y2]
 to show that the   asymotical reduced volume on an asymptotical limit has constant values at all 
times.   Next, in Section 4, we apply the upper bound for $\Delta l$ and the 
convergence lemma in Section 2 to obtain the local $L^2$
strong convergence of $\nabla l$ both in space and in spacetime, and use it to derive the differential inequalities 
[(2.68), Y2] and [(2.69), Y2] in the weak sense for the asymptotical limit. Using this strong convergence we also derive the connection between the two differential inequalitites, 
namely the fact that one becoming an equality forces the other to become an equality.  The remaining steps  are similar to the proof of 
[Theorem 4.9, Y2] in [Y2].

We are grateful to Perelman for providing his lower bound for 
$l$ as presented in [Y2]. We would also like to thank G. Wei 
for  helpful and stimulating discussions. 
 
This paper is based on part of [Y1], whose first verson was posted on the author's webpage 
in February 2004.

\sect{Preliminaries}

Consider a smooth solution $g=g(\tau))$ of the backward Ricci flow
\be  \label{Ricciflow}
\frac{\partial g}{\partial \tau}=2Ric
\ee
on an $n$-dimensional manifold $M$ over an interval $[0, T)$.  We assume that $(M, g(\tau))$ is complete
for each $\tau \in [0, T)$. Choose an arbitary point $p \in M$ as the 
$l$-base, i.e. the reference point for the $l$-function, see [Definition 1, Y1].

We'll follow the notations and conventions in [Y1]. In particular, the volume form of a Riemannian metric 
$h$ on a manifold will be denoted by $dq_h$ or simply $dq$. The distance function of $h$ will 
be denoted by $d_h$ or simply $d$. 



\begin{lem} \label{laplacelemma} Assume that the curvature operator is nonnegative on $[0, T)$. 
Let $\bar \tau \in (0, T)$.   
Then there is a positive constant $C$ depending only on the
dimension $n$ and the magnitude of $\frac{\bar \tau}{T-\bar \tau}$, such that the differential  inequlity $\Delta l \le C\frac{l}{\tau} $ holds true 
for each $\tau \in (0, \bar \tau]$  in 
the weak sense, i.e.
\ba
\int_M \nabla l \cdot \nabla \phi dq \ge -\int_M C \frac{l}{\tau} \phi dq 
\ea
for all nonnegative Lipschitz functions $\phi$ with compact support, as well as 
functions $\phi$ in the Sobolev space $W^{1,2}_{loc}(N)$ with compact support.
(Note that $C$ depends only on $n$ if $T=\infty$.)
\end{lem}
\Pf  By [Theorem 2.20, Y2], the differential inequality  
\ba \label{differential}
\Delta l \le \frac{|\nabla l|^2}{2}+ \frac{R}{2}+\frac{l-n}{2\tau} 
\ea
holds true in the weak sense. By the estimates in [Theorem 2.16, Y2], on $[0, \bar \tau]$  the right hand side 
of (\ref{differential}) is bounded above by $C\frac{l}{\tau}$. The desired conclusion follows. \qed \\

Next we present an analytic lemma regarding strong convergence of 
Sobolev functions. 
 
\begin{lem}  \label{analysislemma1} Let $(N, h)$ be a compact Riemannian manifold, where $h$ denotes the metric. Let $f_k$ be 
a bounded sequence of functions in the Sobolev space $W^{1,2}(N)$ 
such that $f_k-f_{k'} 
\in L^{\infty}(N)$ for all $k, k'$, and that 
\ba \label{gotozero}
f_k-f_{k'} \rightarrow 0
\ea
in $L^{\infty}$ as $k, k' \rightarrow \infty$. 
Furthermore, assume 
that there is a positive constant $C$ such that $\Delta f_k \le C$ 
for all k in the weak sense, i.e. 
\ba \label{super1}
\int_N \nabla f_k \cdot \nabla \phi dq \ge  -C \int \phi dq
\ea
for all nonnegative $\phi \in W^{1,2}(N)$.  Then 
$f_k $ converges strongly in $W^{1,2}(N)$ to a function $f$. 
Consequently, we have 
\ba \label{nablasquare}
\int_N |\nabla f_k|^2 \eta dq \rightarrow \int_N |\nabla f|^2 \eta dq
\ea
for each $\eta \in L^{\infty}(N)$.
The same conclusion holds true if $\Delta f_k \le C$ is replaced 
by $\Delta f_k \ge -C$.
\end{lem}
\Pf  
Since $f_k \in L^2(N)$, (\ref{gotozero}) implies that 
$f_k$ is a Cauchy sequence in $L^2(N)$, and hence converges strongly in 
$L^2(N)$ to a function $f$. Then the boundedness of $f_k$ in $W^{1,2}(N)$ implies that 
$f \in W^{1,2}(N)$ and that $f_k$ converges weakly to $f$. Furthermore, we have 
$f_k-f \in L^{\infty}(N)$ and $f_k-f \rightarrow 0$ in $L^{\infty}(N)$. 
Taking limit in 
(\ref{super1}) we deduce 
\ba \label{super2}
\int_N \nabla f \cdot \nabla \phi dq \ge -C\int_N \phi dq
\ea
for all nonneagtive $\phi \in W^{1,2}(N)$, i.e. $\Delta f \le C$ 
in the weak sense.

We set $\delta_k =\sup |f_k-f|$. Then $\delta_k \rightarrow 0$,
$f-f_k+\delta_k \ge 0$ and $f_k-f+\delta_k \ge 0$. Choosing 
$\phi =f-f_k+\delta_k$ in (\ref{super1}) we deduce
\ba
\int_N \nabla f_k \cdot \nabla f dq -\int_N |\nabla f_k|^2 dq \ge  -C \int (f-f_k+\delta_k) dq.
\ea
Choosing $\phi=f-f_k+\delta_k$ in (\ref{super2}) we deduce
\ba
\int_N |\nabla f|^2 dq -\int_N \nabla f \cdot \nabla f_k dq \ge  -C\int_N (f-f_k+\delta_k) dq.
\ea
It follows that 
\ba \label{bound1}
\int_N |\nabla f|^2 dq -\int_N |\nabla f_k|^2 dq \ge  2C \int(f_k-f-\delta_k) dq.
\ea
On the other hand, choosing $\phi=f_k-f+\delta_k$ in (\ref{super1}) and 
(\ref{super2}) we infer 
\ba \label{bound2}
\int_N |\nabla f|^2 dq -\int_N |\nabla f_k|^2 dq \le 2C \int(f_k-f+\delta_k) dq.
\ea
Obviously, (\ref{bound1}) and (\ref{bound2}) imply 
\ba \label{nablasquare1}
\int_N |\nabla f_k|^2  dq \rightarrow \int_N |\nabla f|^2 dq.
\ea
It is easy to see that (\ref{nablasquare1}) implies that 
$f_k$ converges strongly to $f$ in $W^{1,2}(N)$. Indeed, we have 
\ba
\int_N |\nabla f_k-\nabla f|^2 dq =\int_N |\nabla f_k|^2 dq 
+\int_N |\nabla f|^2 dq -2\int_N \nabla f_k \cdot \nabla f dq.
\ea
The right hand side is easily seen to converge to zero.

Now we derive (\ref{nablasquare}). We have 
\ba 
\int_N |\nabla f_k|^2 \eta dq-\int_N |\nabla f|^2 \eta dq &=&
\int_N (\nabla f_k-\nabla f) \cdot \nabla f \eta dq + \nonumber \\&&
\int_N (\nabla f_k-\nabla f) \cdot \nabla f_k \eta dq.
\ea
Applying Cauchy-Schwarz inequality and the strong convergence of 
$\nabla f_k$ to $\nabla f$ in $L^2$ we then arrive at (\ref{nablasquare}).

The proof in the case $\Delta f_k \ge -C$ is similar.
\qed \\

It is easy to extend this lemma to the situation of 
noncompact manifolds. For convenience, we formulate the result 
in the set-up of pointed manifolds.

\begin{lem} \label{analysislemma2} Let $(N, h, p)$ be a pointed complete noncompact Riemannian 
manifold. Let $f_k$ be a locally bounded sequence of functions in $W^{1,2}_{loc}(N)$ 
such that 
the $\Delta f_k$ are locally unformly bounded from above in the weak sense, namely  
 there exist for each $A>0$ positive constants $C_1(A)$ and $C_2(A)$ such that for each $k$
 the $W^{1,2}$-norm of $f_k$ on the geodesic ball $d(p, \cdot)<A$ is bounded above by $C_1(A)$, and 
$\Delta f_k \le C_{A}$ holds true on $d(p, \cdot)<A$ in the weak sense, i.e. 
\ba \label{ksuper1}
\int_{N_k} \nabla f_k \cdot \nabla \phi dq \ge -C_A \int_{N} \phi dq
\ea
for all nonnegative $\phi \in W^{1,2}_{loc}(N)$ with support 
contained in $d(p, \cdot)<A$.  Furthermore, assume that 
$f_k-f_{k'} \in L^{\infty}_{loc}(N)$ for all $k, k'$, and 
$f_k-f_{k'} \rightarrow 0$ in $L^{\infty}_{loc}(N)$. 
Then $f_k$ converges strongly in $W^{1,2}_{loc}$ to a function $f$. Consequently, 
we have
\ba \label{nablasquare2}
\int_{N} |\nabla f_k|^2 \eta dq \rightarrow \int_{N}|\nabla f|^2 \eta dq
\ea
for each $\eta \in L^{\infty}(N)$ with compact support. 

The same conclusion holds true if the $\Delta f_k$ are assumed to be 
locally uniformly bounded from below in the weak sense instead of from above. 
\end{lem}

\Pf 
First, $f_k$ converges to some $f$ weakly in $W^{1,2}_{loc}(N)$. Moreover, 
$f_k-f$ converges to zero in $L^{\infty}_{loc}(N)$.  As in the proof of Lemma \ref{analysislemma1}
we obtain by taking limit in  (\ref{ksuper1}) 
\ba \label{ksuper2}
\int_{N} \nabla f \cdot \nabla \phi dq \ge -C_A \int_{N} \phi dq
\ea
for all $A>0$ and all nonnegative $\phi \in W^{1,2}_{loc}(N)$ with support contained in 
$d(p, \cdot)<A$. For a fixed 
$A>0$, we set $\delta_k=\sup_{d(p, \cdot)\le 2A} |f_k-f|$. 
Let $\psi$ be a nonnegative Lipschitz function on $N$ which takes the value 
$1$ on $d(p, \cdot) \le A$ and is zero on $d(p, \cdot) \ge 2A$. 
Choosing $\phi=\psi(f-f_k+\delta_k)$ in (\ref{ksuper1}) and (\ref{ksuper2}) we deduce 
\ba \label{newbound1}
\int_N |\nabla f|^2 dq -\int_N |\nabla f_k|^2 dq  &\ge& 2C_{2A} \int \psi (f_k-f-\delta_k) dq
+\int_N \nabla f_k \cdot \nabla \psi (f_k-f-\delta_k) dq   \nonumber \\
&& +\int_N \nabla f \cdot \nabla \psi (f_k-f-\delta_k) dq.
\ea
On the other hand, choosing $\phi=\psi(f_k-f+\delta_k)$ in (\ref{ksuper1}) and 
(\ref{ksuper2}) we infer
\ba \label{newbound2}
&& \int_N |\nabla f|^2\psi  dq - \int_N |\nabla f_k|^2 \psi dq \le 2C_{2A} \int_N \psi (f-f_k+\delta_k) dq
\nonumber \\ &+& \int_N \nabla f_k \cdot \nabla \psi (f_k-f+\delta_k) dq 
+ \int_N \nabla f \cdot \nabla \psi (f_k-f+\delta_k) dq.
\ea
We can use the local boundedness of $f_k$ in $W^{1,2}_{loc}$ and Cauchy-Schwarz inequality to handle the second and third  terms on the right
hand side of (\ref{newbound1}) and (\ref{newbound2}). It follows that 
$\int_N|\nabla f_k|^2 \psi dq \rightarrow \int_N |\nabla f|^2 \psi dq$. Then
\ba 
\int_{d(p, \cdot) \le A} |\nabla f_k-\nabla f|^2 dq &\le& \int_N |\nabla f_k-\nabla f|^2 \psi dq 
=\int_N |\nabla f_k|^2 \psi dq +\int_N |\nabla f|^2 \psi dq \nonumber \\
&& - 2\int_N \nabla f_k \cdot \nabla f \psi dq \rightarrow 0.
\ea

The derivation of (\ref{nablasquare2}) is similar to that of (\ref{nablasquare}) in the 
proof of Lemma \ref{analysislemma1}.

The case that $\Delta f_k$ are locally uniformly bounded from below in the weak sense is similar.
\qed
\\

Next we extend the above results to the situation of a sequence of 
converging Riemannian manifolds.  For simplicity of formulation 
we state a weaker version of the result. 

\begin{lem} \label{analysislemma3} Let $(N_k, h_k, p_k)$ be a sequence of pointed complete Riemannian 
manifolds 
converging in $C^1$ to a pointed complete Riemannian manifold 
$(N, h, p)$. Let $f_k \in W^{1,2}_{loc}(N_k)\cap L^{\infty}_{loc}(N_k)$ which converges in $L^{\infty}_{loc}$ and 
weakly in $W^{1,2}_{loc}$ to a function $f \in W^{1,2}_{loc}(N) \cap L^{\infty}_{loc}(N)$.
Furthermore, we assume that the $\Delta f_k$ are locally uniformly bounded from above in the weak sense or 
locally uniformly bounded from below in the weak sense. 
Then $f_k$ converges strongly in $W^{1,2}_{loc}$ to $f$. 
Consequently, 
we have
\ba
\int_{N_k}|\nabla f_k|^2 \eta_k dq \rightarrow \int_{N}|\nabla f|^2 \eta dq,
\ea
whenever $\eta \in L^{\infty}(N)$ and $\eta_k \in L^{\infty}(N_k)$ with support contained in 
$d_{h_k}(p_k, \cdot) <A$ for a fixed $A>0$, such that $\eta_k$ converges 
to $\eta$ in $L^{\infty}$.
\end{lem}

\Pf The proof of Lemma \ref{analysislemma2} can easily be carried over. Note that the pointed convergence of 
$(N_k, h_k, p_k)$ to $(N, h, p)$ involves suitable embeddings from increasing domains of $N_k$ into 
$N$. Using these embeddings the functions $f_k$ can be transplanted to increasing domains of $N$. Then 
the difference between $f_k$ and $f$ can be measured.
\qed

\sect{Asymptotic Limits of $\kappa$-Solutions I}

Let $\tilde g(t), -\infty<t\leq 0$ be a $\kappa$-solution for some $\kappa>0$ on
a manifold $M$.
Recall [P1] that this means that $\tilde g(t)$ is a smooth solution of the Ricci flow
\be
\frac{\partial g}{\partial t}=-2Ric
\ee
on $M \times (-\infty, 0]$, such that for each $t$, the metric $\tilde g(t)$ is
a complete non-flat metric of bounded and  nonnegative curvature
operator. Moreover,
each $\tilde g(t)$ is $\kappa$-noncollapsed on all scales. (See [P1] for the 
definition of $\kappa$-noncollapsedness.)  To understand the
structures of $\tilde g(t)$, we analyse its rescaled asymptotical limits at the time
infinity.  One needs to use a blow-down rescaling, because the Ricci flow
equation and nonnegative curvature imply that $\tilde g(t)$ expands as $t$
decreases.

Pick an arbitrary time $t_0\le 0$.
Then $g(\tau)=\tilde g(t_0-\tau), \tau \in [0, \infty)$ is
a solution of the backward Ricci flow (\ref{Ricciflow}) on $M \times [0, \infty)$. Obviously, 
$g$ retains the properties of $\tilde g$, namely $g$ is $\kappa$-noncollapsed
on all scales and has bounded and nonnegatuve curvature operator.



Next we
choose an arbitary $p \in M$ as the $l$-base, i.e. the reference point for 
the $l$-function $l$ and the reduced volume $\tilde V$ (see [Definition 1, Y2]). We set
\ba
\tilde V^{\infty}=\lim_{\tau \rightarrow \infty} \tilde V(\tau).
\ea
By [Theorem 4.5, Y2], this limits exists and is finite.

For $a>0$, we set $g_{a}(\tau)=\frac{1}{a}g(a \tau)$. 
Choose a base point $\bar q \in M$ we then have a pointed  
solution of the backward Ricci flow
$(g_a, M \times (0, \infty), \bar q)$. The following
theorem is precisely [Proposition 11.2, P1].



\begin{theo} \label{theorem1}  Let $\tau_k \rightarrow \infty$
be given. For each $\tau_k$, let $x(\tau_k)$ be a minimum point of 
$l(\cdot, \tau)$. Then the pointed flows $(g_{\tau_k}, M \times (0, \infty),
x(\tau_k))$  subconverge smoothly to pointed nonflat gradient shrinking solitons
$(g_{\infty}, M_{\infty} \times (0, \infty), x_{\infty})$.
\end{theo}




In the remainder of this section, we establish the convergence part of 
Theorem \ref{theorem1}, and 
obtain some basic properties for the limit solutions. For convenience, we formulate the  
convergence part as a proposition.

\begin{prop} The pointed flows $(g_{\tau_k}, M \times (0, \infty), x(\tau_k))$ 
subconverge 
smoothly to pointed solutions $(g_{\infty}, M_{\infty} \times (0, \infty), x_{\infty})$
of the backward Ricci flow.
The limits will be called asymptotical limits of $\tilde g$ or $g$.
\end{prop}



\Pf 
 The various quantities associated with $g_{\tau_k}$ will be 
indicated by the subscript $k$ or $g_{\tau_k}$, e.g. $l_k=l_{g_{\tau_k}}$, $\tilde V_k=\tilde
V_{g_{\tau_k}}$ (the $l$-function and reduced volume associated with
$g_{\tau_k}$), and $d_k=d_{g_{\tau_k}}$.

First note that by [Lemma 2.3, Y2] or [Lemma 3.2, Y2] the minimum points $x(\tau_k)$ exist. 
By [Lemma 3.1, Y2] and the scaling invariance of the $l$-function we 
have 
\ba \label{dan}
l_k(x(\tau_k), 1) \le \frac{n}{2}.
\ea
By this estimate and [Lemma 3.2, Y2] we infer  
\be  \label{rescalelestimateI}
l_k (q, 1)\leq C_2 d^2_{k}(x(\tau_k), q, 1)+n
\ee
for all $q \in M$, where $C_2$ is a positive constant depending only on the dimension $n$.
Then it follows from the 
Harnack inequality [(2.52), Y2]   that 
\be  \label{rescalelestimateII}
l_k(q, \tau)\leq \tau^{\pm C}
(C_2(d^2_{k}(x(\tau_k), q, 1)+n),
\ee
where $\pm=+$ if $\tau\ge 1$, 
$\pm=-$ if $\tau<1$, 
and 
$C$ is a positive constant depending only on $n$.
Consequently, we obtain from [Theorem 2.16, Y2]  and the nonnegativity of curvature
operator
the estimate
\be \label{curvature}
|Rm|_k(q, \tau) \leq C \tau^{-1\pm C}
(d^2_{k}(x(\tau_k), q, 1)+1).
\ee

By the $\kappa$-noncollapsing assumption and (\ref{curvature}), we obtain
for each $k$ and $\tau \in (0, \infty)$ local injectivity radius estimates for $g_{\tau_k}(\tau)$
which depend only on $d_{k}(x(\tau_k),
\cdot, 1)$, $\tau$ (in terms of a positive lower bound and an upper bound 
of $\tau$), $\kappa$ and $n$. It is then straightforward to apply Gromov-Cheeger-Hamilton compactness [H]
to obtain pointed smooth subconvergence of $(g_{\tau_k}, M \times (0, \infty), x(\tau_k))$
to pointed solutions $(g_{\infty}, M_{\infty} \times (0, \infty), x_{\infty})$ of the 
backward Ricci flow.
\qed 
\\

Let $(g_{\infty}, M_{\infty} \times (0, \infty), x_{\infty})$ be an asymptotical limit. 
The corresponding converging subsequence
will still be denoted by $g_{\tau_k}$.  By [Theorem 2.16, Y2], the recaling invariance
 and (\ref{rescalelestimateII}), the $l_{k}$ subconverge locally uniformly to locally Lipschitz functions
$l_{\infty}$ on $M_{\infty} \times (0, \infty)$. 
\\

\noindent {\bf Definition 1} The limits $l_{\infty}$ will be called 
{\it limit $l$-functions}. 
\\

\begin{lem} Let $(g_{\infty}, M_{\infty} \times (0, \infty), x_{\infty})$ be 
an asymptotical limit and $l_{\infty}$ a limit $l$-function. Then  
$l_{\infty}(\cdot, \tau)$ is locally Lipschitz 
for every $\tau$, and $l_{\infty}(q, \cdot)$ is locally 
Lipschitz for every $q$.  Moreover, we have 
\ba \label{limitboundA}
R_{{\infty}}
\le C\frac{l_{\infty}}{\tau}
\ea
everywhere, 
\ba \label{limitboundB}
|\nabla  l_{\infty}|
\le C\frac{l_{\infty}}{\tau}
\ea
almost everywhere in $M_{\infty}$ for each $\tau \in (0, \infty)$ (the $\nabla$ and the norm are those 
of  $g_{\infty}$), 
\ba \label{limitboundC}
|\frac{\partial l_{\infty}}{\partial \tau}|
\le C\frac{l_{\infty}}{\tau}
\ea
almost everywhere in $(0, \infty)$ for each $q \in M_{\infty}$, and 
\ba \label{limitboundD}
(\frac{\tau_1}{\tau_2})^C \le \frac{l_{\infty}(q, \tau_2)}{l_{\infty}(q, \tau_1)}\leq
(\frac{\tau_2}{\tau_1})^C
\ea
for all $q \in M$ and $0<\tau_1<\tau_2$,
where $R_{\infty}=R_{g_{\infty}}$ (the scalar curvature of $g_{\infty}$)
and  $C$ depends only on $n$. 
\end{lem}
\Pf These follow from the local uniform convergence of 
$l_k$ to $l_{\infty}$ and the corresponding properties of the $l$-function 
as given in [Theorem 2.16, Y2]. The estimates for the gradient and the $\tau$-derivative 
are obtained in terms of estimating the relevant Lipschitz constants.
\qed

\begin{lem} \label{l-distance}
We have 
\be \label{onecenter}
C_1 \frac{d_k^2(x(\tau_k),q, \tau)}{\tau}-\frac{n}{2}\tau^{\pm C}-1 
\le l_k(q, \tau) \le C_2 \frac{d_k^2(x(\tau_k), q, \tau)}{\tau}+n\tau^{\pm C}  
\ee
for all $q \in M, \tau \in (0, \infty)$ and 
\be  \label{l-estimate}
C_1\frac{d^2_{\infty}(x_{\infty}, q, \tau)}{\tau}-\frac{n}{2}\tau^{\pm C}-1 \le l_{\infty} (q, \tau)\leq C_2 \frac{d^2_{\infty}(x_{\infty}, q, \tau)}{\tau}+n \tau^{\pm C}
\ee
for all $q \in M_{\infty}, \tau \in (0, \infty)$, where $d_{\infty}$ denotes the distance 
function of $g_{\infty}$, and the positive constant $C_1$, $C_2$ and 
$C$ depend only on the dimension $n$.
\end{lem}
\Pf  By the Harnack inequality [(2.52), Y2] and (\ref{dan}) we deduce
\ba
l_k(x(\tau_k), \tau) \le \frac{n}{2}\tau^{\pm C}.
\ea
Applying [Lemma 3.2, Y2] we then arrive at (\ref{onecenter}). Taking limit 
leads to (\ref{l-estimate}). (Note that a similar estimate was employed in 
the proof of [Theorem 3.3, Y2] in [Y2].) \qed \\

Henceforth, we'll fix a limit $l$-function $l_{\infty}$ for each 
asymptotical limit. Obviously, the conclusions we 
derive are valid for each choice of the limit $l$-function.\\

\noindent {\bf Definition 2} Let $(g_{\infty}, M_{\infty} \times (0, \infty), x_{\infty})$
be an asymptotical limit together with a limit $l$-function $l_{\infty}$. We define the {\it asymptotical reduced volume}
to be  
\ba
\tilde V_{\infty}(\tau)=\int_{M_{\infty}} e^{-l_{\infty}} \tau^{-\frac{n}{2}}
dq,
\ea
where $dq=dq_{g_{\infty}}$. 
\\

\begin{lem} \label{volumelemma} 
$\tilde V_k(\tau)$ converges to $\tilde V_{\infty}(\tau)$ for each $\tau>0$
and hence there holds
\ba \label{volumevalue}
\tilde V_{\infty}(\tau)=\tilde V^{\infty}
\ea
for each $\tau$. 
\end{lem}

\Pf 
Fix $\tau>0$.  For each $A>0$ we have 
\ba \label{split1}
\tilde V_k(\tau)=\int_{d_k(x(\tau_k), \cdot, \tau) \le A} e^{-l_k} \tau^{-\frac{n}{2}}dq_k
+\int_{d_k(x(\tau_k), \cdot, \tau) > A} e^{-l_k} \tau^{-\frac{n}{2}}dq_k
\ea
and 
\ba \label{split2}
\tilde V_{\infty}(\tau)=\int_{d_{\infty}(x_{\infty}, \cdot, \tau) \le A} e^{-l_{\infty}} \tau^{-\frac{n}{2}}dq
+\int_{d_{\infty}(x_{\infty}, \cdot, \tau) > A} e^{-l_{\infty}} \tau^{-\frac{n}{2}}dq.
\ea
As in the proof of [Theorem 3.3, Y2] in [Y2], the lower bound in (\ref{onecenter})
and volume comparison imply 
\ba \label{smallI}
\int_{d_k(x(\tau_k), \cdot, \tau) > A} e^{-l_k} \tau^{-\frac{n}{2}}dq_k \le Ce^{-cA}
\ea
for suitable positive constants $C$ and $c$ depending only on $n$.  Similarly, the lower bound in (\ref{l-estimate}) 
and volume comparison imply
\ba \label{smallII}
\int_{d_{\infty}(x_{\infty}, \cdot, \tau) > A} e^{-l_{\infty}} \tau^{-\frac{n}{2}}dq \le Ce^{-cA}.
\ea
(We can choose the constants here to be the same as in (\ref{smallI}).) 
By the pointed convergence of $(g_{\tau_k}, M \times (0, \infty), x(\tau_k))$ to $(g_{\infty},M_{\infty} \times
(0, \infty), x_{\infty})$
and the convergence of $l_k$ to $l_{\infty}$ we then infer 
$\lim |\tilde V_k(\tau)-\tilde V_{\infty}(\tau)| \le 2Ce^{-cA}$. Since $A$ is arbitary, we conclude 
that $\lim \tilde V_k(\tau)=\tilde V_{\infty}(\tau)$.  But $\tilde V_k(\tau)=\tilde V(\tau_k \tau) 
\rightarrow \tilde V^{\infty}$. Hence we arrive at (\ref{volumevalue}). \qed

\sect{Asymtotical Limits of $\kappa$-solutions II}

The goal of this section is to establish the soliton characterization part 
of Theorem \ref{theorem1}. We formulate it
as a proposition.  See [Y2] for the terminologies regarding gradient 
shrinking solitons used below. 

\begin{prop} Each asymptotical limit $(g_{\infty}, M_{\infty} \times (0, \infty), x_{\infty})$ 
is a pointed nonflat gradient shrinking soliton with time origin $0$. Moreover, 
the limit $l$-functions $l_{\infty}$ are potential functions.
\end{prop}

To prove this proposition, we need a few lemmas. We fix an asymptotical 
limit $(g_{\infty}, M_{\infty} \times (0, \infty), x_{\infty})$ together with a
limit $l$-function $l_{\infty}$.  We have the corresponding converging 
subsequences $(g_{\tau_k}, M \times (0, \infty), x(\tau_k))$ and $l_k$.

\begin{lem} \label{strong1} For each $\tau>0$, $l_k(\cdot, \tau)$ 
converges strongly to $l_{\infty}(\cdot, \tau)$ in $W^{1,2}_{loc}$. 
Consequently, we have
\ba \label{convergenceA} 
\int_{M, g_{\tau_k}} |\nabla l_k|^2 \phi_k dq|_{\tau} \rightarrow \int_{M_{\infty}} |\nabla l_{\infty}|^2 \phi dq,
\ea
whenever $\phi \in L^{\infty}(M_{\infty})$ and $\phi_k \in L^{\infty}(M)$ with support 
contained in $d_k(x(\tau_k), \cdot, \tau) <A$ for a fixed $A>0$, such that 
$\phi_k$ converges in $L^{\infty}$ to $\phi$. 
\end{lem}
\Pf  By [Theorem 2.16, Y2] and (\ref{rescalelestimateII}), 
$l_k$ is a locally bounded sequence in $W^{1,2}_{loc}$. Since it 
converges to $l_{\infty}$ in $L^{\infty}_{loc}$, it also converges 
to $l_{\infty}$ weakly in $W^{1,2}_{loc}$.   Hence the conclusions of the lemma  follow from Lemma \ref{laplacelemma}, Lemma \ref{analysislemma3} and Lemma \ref{l-distance}. 
\qed \\

We also need the $L^2_{loc}$ strong convergence of $\nabla l_k$ over the 
spacetime $M \times (0, \infty)$. 

\begin{lem} \label{strong2} The vector fields  $\nabla l_k$ on $M \times (0, \infty)$ converge  to 
the vector field $\nabla l_{\infty}$ on $M_{\infty} \times (0, \infty)$ strongly in $L^2_{loc}$. Consequently,
we have for arbitary $\tau_2>\tau_1>0$ 
\ba \label{convergenceB} 
\int_{\tau_1}^{\tau_2}\int_{M, g_{\tau_k}} |\nabla l_k|^2 \phi_k dqd\tau \rightarrow 
\int_{\tau_1}^{\tau_2} \int_{M_{\infty}} |\nabla l_{\infty}|^2 \phi dqd\tau,
\ea
whenever $\phi \in L^{\infty}(M_{\infty}\times [\tau_1, \tau_2])$ and $\phi_k 
\in L^{\infty}(M \times [\tau_1, \tau_2])$ with support contained in 
$\{d_k(x(\tau_k), \cdot, 1) < A \} \times [\tau_1, \tau_2]$ for a fixed $A>0$, such that 
$\phi_k$ converges to $\phi$ in $L^{\infty}$.
\end{lem}
\Pf By [Theorem 2.16, Y2] and (\ref{rescalelestimateII}), $l_k$ is a locally bounded sequence in $W^{1,2}_{loc}$ over 
the spacetime.  It follows that $l_k$ converges weakly in $W^{1,2}_{loc}$ over the 
spacetime. Now we can argue as in the proof of Lemma \ref{analysislemma2}, employing 
integrations over the spacetime with respect to the volume form $dqd\tau$. 
 Note that the involved $\Delta$ and $\nabla$ are along the space with a time 
 dependence on $\tau$.  \qed

\begin{lem} \label{equations} The equation 
\ba \label{theequation}
\frac{\partial l_\infty}{\partial \tau}-\frac{R_{\infty}}{2}+\frac{|\nabla l_{\infty}|^2}{2}+
\frac{l_{\infty}}{2\tau}=0
\ea
holds true almost everywhere on $M_{\infty} \times (0, \infty)$. The inequality
\ba
\Delta l_{\infty}-\frac{|\nabla l_{\infty}|^2}{2}+\frac{R_{\infty}}{2}+\frac{l_{\infty}-n}{2\tau}\le 0
\ea
holds true for each $\tau>0$ in the weak sense, i.e. 
\ba \label{thelaplace}
\int_{M_{\infty}} \{-\nabla l_{\infty} \cdot \nabla \phi  
+\frac{1}{2}(-|\nabla l_{\infty}|^2+R_{\infty}+\frac{l_{\infty}-n}{\tau})\phi \} dq
\le 0
\ea
for all nonnegative Lipschitz functions 
$\phi$ with compact support.  Finally, the inequality 
\ba  \frac{\partial l_{\infty}}{\partial \tau}-\Delta l_{\infty}+|\nabla l_{\infty}|^2-R_{\infty}+\frac{n}{2\tau} \ge 0
\ea
holds true on $M_{\infty} \times (0, \infty)$ when $\Delta$ is interpreted in the weak sense, i.e. 
\be  \label{Qnonnegative}
Q_{\tau_1, \tau_2}(\phi) \geq 0 \ee for arbitray $\tau_2>\tau_1>0$
and nonnegative Lipschitz functions $\phi$ on $M_{\infty} \times
[\tau_1, \tau_2]$  with compact support, where \be Q_{\tau_1,
\tau_2}(\phi)= \int_{\tau_1}^{\tau_2}\int_{M_{\infty}}
\{\nabla l_{\infty} \cdot  \nabla \phi
+(\frac{\partial l_{\infty}}{\partial \tau}+ |\nabla
l_{\infty}|^2-R_{{\infty}}+\frac{n}{2\tau})\phi \} dq d\tau. \ee
\end{lem}
\Pf First consider a fixed $\tau>0$. Let $\phi$ be a nonnegative smooth function on $M_{\infty}$ with compact support. 
Let $\phi_k$ be a sequence of nonnegative smooth functions on $M$ such that the support of 
$\phi_k$ is contained in $d_k(x(\tau_k), \cdot, \tau)<A$ for a fixed $A>0$, and that 
$\phi_k$ converges smoothly to $\phi$.  By [Theorem 2.20, Y2] we have for each $k$
\ba \label{kequation}
\int_{M, g_{\tau_k}} \{-\nabla l_k \cdot \nabla \phi_k +(-\frac{|\nabla l_k|^2}{2}+\frac{R_k}{2}+\frac{l_k-n}{2\tau})\phi_k\}
dq\leq 0
\ea
 On account of Lemma \ref{strong1} and the 
pointed convergence of $(g_{\tau_k}, M\times (0, \infty), x(\tau_k))$ to $(g_{\infty},
M_{\infty} \times (0, \infty), x_{\infty})$ we can take limit in  (\ref{kequation}) to 
arrive at (\ref{thelaplace}) with the above $\phi$.  Then it follows that 
(\ref{thelaplace}) also holds true for nonnegative Lipschitz or Sobolev functions
$\phi$ with compact support. 

Next consider $\tau_2>\tau_1>0$
and a smooth nonnegative function $\phi$ on $M_{\infty} \times [\tau_1, \tau_2]$ with 
compact support.  Let $\phi_k$ be a sequence of smooth functions on $M \times [\tau_1, \tau_2]$ 
such that the support of $\phi_k$ is contained in $d_k(x(\tau_k), \cdot, 1) <A$ for 
a fixed $A>0$, and that $\phi_k$ converges smoothly to $\phi$. 
By [Theorem 2.20, Y2] we have
\ba \label{kheat}
\int_{\tau_1}^{\tau_2} \int_{M, g_{\tau_k}} \{\nabla l_k \cdot
\nabla \phi_k+(\frac{\partial l_{k}}{\partial \tau}+|\nabla l_k|^2-R_k+\frac{n}{2\tau})\phi_k
\} dq d\tau \geq 0.
\ea
Since $\frac{\partial}{\partial \tau} dq=R_kdq$ we can rewrite this equation as follows
\ba \label{rewrite}
&& \int_{M, g_{\tau_k}} l_k  \phi_k dq|_{\tau_2}-\int_{M, g_{\tau_k}} l_k \phi_k dq|_{\tau_1}+  \nonumber \\
&&\int_{\tau_1}^{\tau_2} \int_{M, g_{\tau_k}} \{\nabla l_k \cdot
\nabla \phi_k+(-l_k \frac{\partial \phi_{k}}{\partial \tau}+|\nabla l_k|^2-2R_k+\frac{n}{2\tau})\phi_k
\} dq d\tau \geq 0.
\ea
On account of Lemma \ref{strong2}, the local uniform convergence of $l_k$ to 
$l_{\infty}$ and the pointed convergence of $(g_{\tau_k}, M\times (0, \infty), x(\tau_k))$
we can take limit in (\ref{rewrite}) to arrive at 
\ba
&&\int_{M_{\infty}} l_{\infty}  \phi dq|_{\tau_2}-\int_{M_{\infty}} l_{\infty} \phi dq|_{\tau_1}+
\nonumber \\ && \int_{\tau_1}^{\tau_2} \int_{M_{\infty}} \{\nabla l_{\infty} \cdot
\nabla \phi+(-l_{\infty} \frac{\partial \phi}{\partial \tau}+|\nabla l_{\infty}|^2-2R_{\infty}+\frac{n}{2\tau})\phi
\} dq d\tau \geq 0,
\ea
which is equivalent to (\ref{Qnonnegative}) for the above $\phi$. Then (\ref{Qnonnegative}) also 
holds true for nonnegative Lipschitz or Sobolev $\phi$ with compact support.

By [Lemma 2.19, Y2] and [Lemma 2.15, Y2] we have 
\ba \label{beforeequation}
\frac{\partial l_k}{\partial \tau}-\frac{R_{k}}{2}+\frac{|\nabla l_{k}|^2}{2}+
\frac{l_{k}}{2\tau}=0
\ea 
almost everywhere on $M \times (0, \infty)$. We turn this equation into an integral equation 
in the same fashion as (\ref{rewrite}). Taking limit we then arrive at a corresponding integral 
equation for $l_{\infty}$, which implies that  (\ref{theequation}) holds true almost 
everywhere. \qed \\


\noindent {\it Proof of Proposition 4.1} \\

Fix an asymptotical limit $(g_{\infty}, M_{\infty} \times (0, \infty), x_{\infty})$ 
together with a limit $l$-function $l_{\infty}$ as above. By (\ref{Qnonnegative}), Lemma \ref{l-distance} and 
the arguments in the proof of [Theorem 3.3, Y2] in [Y2]
we infer that $Q_{\tau_1, \tau_2}(\phi) \ge 0$ holds true for nonnegative admissible locally Lipschitz $\phi$ 
on $M_{\infty} \times [\tau_1, \tau_2]$, where ``admissible" means 
that $ \phi \le C \tau^{-\frac{n}{2}}
e^{-l_{\infty}}$ and $|\nabla \phi| \le C \tau^{-\frac{n+1}{2}}\sqrt{l_{\infty}} e^{-l_{\infty}}$ for all $q \in M_{\infty}$ and $\tau
\in [\tau_1, \tau_2]$ and some positive constant $C$
depending on $\phi$. On the other hand, we can argue as in the proof of [Lemma 4.6, Y2] in [Y2] 
to deduce 
\ba
\tilde V_{\infty}(\tau_2)-\tilde V_{\infty}(\tau_1)=
-\int_{\tau_1}^{\tau_2}\int_{M_{\infty}} (\frac{\partial l_{\infty}}{\partial \tau}
-R_{\infty}+\frac{n}{2\tau})e^{-l_{\infty}}\tau^{-\frac{n}{2}} dqd\tau. 
\ea
By Lemma \ref{volumelemma} we then infer that $Q(\tau^{-\frac{n}{2}}e^{-l_{\infty}})=0$.
As in the proof of [Theorem 4.9, Y2] in [Y2], we then deduce that $Q_{\tau_1, \tau_2}(\phi)=0$
for all admissible locally Lipschitz $\phi$.  By basic regularity 
theory for parabolic equations we derive that $l_{\infty}$ is smooth  and 
the equation 
\ba \label{newheat}
\frac{\partial l_{\infty}}{\partial \tau}-\Delta l_{\infty}+|\nabla l_{\infty}|^2-R_{\infty}+\frac{n}{2\tau} = 0
\ea
holds true eveywhere. By Lemma \ref{equations}, the equation (\ref{theequation}) also holds 
true everywhere.  Taking the difference between  
(\ref{theequation}) and (\ref{newheat}) we then arrive at the equation
\ba \label{newlaplace}
\Delta l_{\infty}-\frac{|\nabla l_{\infty}|^2}{2}+\frac{R_{\infty}}{2}+\frac{l_{\infty}-n}{2\tau}= 0.
\ea
By (\ref{newheat}), (\ref{newlaplace}) and [Lemma 4.8,  Y2] we finally conclude that 
$g_{\infty}$ is a gradient shrinking soliton with time origin $0$ and potential function $l_{\infty}$.

Now we show that  $g_{\infty}$  is
nonflat.  Assume that $g_{\infty}$ is flat.  Then the soliton
equation [(4.16), Y2] implies \be \label{flat1} \nabla^2_{g_{\infty}}
l_{\infty}= \frac{1}{2\tau}g_{\infty}. \ee It follows that
$l_{\infty}(\cdot, 1)$ is a strictly convex function with a unique
minimum point (by the construction it has a minimum point
$x_{\infty}$). By Lemma \ref{l-distance}, its level sets are
compact. Employing the gradient flow of $l_{\infty}(\cdot, 1)$ 
one readily shows that  $M_{\infty}$ is diffeomorphic to ${\bf R}^n$ and
then each $(M_{\infty}, g_{\infty}(\tau))$ is isometric to ${\bf
R}^n$. Next we observe that (\ref{flat1}) implies that
$\Delta_{g_{\infty}} l_{\infty}=\frac{n}{2\tau}$, which together
with (\ref{newlaplace}) leads to \be \label{flat2} |\nabla
l_{\infty}|^2=\frac{l_{\infty}}{\tau}. \ee For a fixed $\tau$, we
identify $(M_{\infty}, g_{\infty})$ with ${\bf R}^n$ via an
isometry. Let $x$ denote the coordinates on ${\bf R}^n$. Then
(\ref{flat1}) implies that \be  \label{flat3} \nabla l_{\infty}(x,
\tau)=\nabla \frac{|x|^2}{4\tau}+v(\tau) \ee for a constant vector
$v(\tau)$.  Hence \be  \label{flat4} l_{\infty}(x,
\tau)=\frac{|x|^2}{4\tau}+v(\tau) \cdot x+c(\tau) \ee for a
constant $c(\tau)$. A simple calculation using (\ref{flat2}),
(\ref{flat3}) and (\ref{flat4}) yields $c(\tau)=\tau |v|^2$. It
follows that $l_{\infty}(x, \tau)=\frac{|x+2\tau v|^2}{4\tau}$ and
hence $\tilde V_{{\infty}}(\tau)=(4\pi)^{\frac{n}{2}}$ for each $\tau$. But $\tilde
V_{{\infty}}(\tau)=\tilde V^{\infty}<(4\pi)^{\frac{n}{2}}$ by [Theorem 4.4,
Y2] and [Theorem 4.5, Y2]. Thus we arrive at a contradiction. 
\qed

Department of Mathematics,
 University of California,
Santa Barbara, CA 93106

yer@math.ucsb.edu


\begin{thebibliography}{DS}
\bibitem[H]{h1}  R.~Hamilton,
{\em A compactness property for solutions of the Ricci flow,}
Amer.~J.~Math.~117 (1995), 545-572.
\bibitem[P1]{p1} G.~Perelman,
{\em The entropy formula for the Ricci flow and its geometric
applications},
arXiv:math.DG/0211159.
\bibitem[P2]{p2} G.~Perelman, {\em Ricci flow with surgeries on three-manifolds},
arXiv:math.DG/0303109.
\bibitem[Y1]{yy1} R.~Ye,
{\em On the $l$-function and the reduced volume of Perelman.}
\bibitem[Y2]{yy2} R.~Ye,
{\em On the $l$-function and the reduced volume of Perelman I.}
\end{thebibliography}
\end{document}